\documentstyle[12pt]{article}

  \newcommand{\const}{\rm const}

\begin{document}

 \begin{center}

 \ {\bf Weighted  Grand Lebesgue Spaces norm estimation}\par

\vspace{4mm}

{\bf  for Hardy - Sobolev - Poincare - Wirtinger operators.} \par

\vspace{5mm}

{\bf  M.R.Formica, E.Ostrovsky, and L.Sirota. } \par

 \vspace{5mm}

 \end{center}

 \ Universit\`{a} degli Studi di Napoli Parthenope, via Generale Parisi 13, Palazzo Pacanowsky, 80132,
Napoli, Italy. \\

e-mail: mara.formica@uniparthenope.it \\

\vspace{3mm}

Department of Mathematics and Statistics, Bar-Ilan University,
59200, Ramat Gan, Israel. \\

\vspace{4mm}

e-mail: eugostrovsky@list.ru\\
Department of Mathematics and Statistics, Bar-Ilan University,\\
59200, Ramat Gan, Israel.

\vspace{4mm}

e-mail: sirota3@bezeqint.net \\

\vspace{5mm}

\begin{center}

 \ {\it Abstract.}

 \end{center}

 \ We  extend the classical Hardy - Sobolev - Poincare - Wirtinger inequalities from the ordinary Lebesgue - Riesz
 spaces into the Grand Lebesgue ones, with exact constants evaluation.\par

 \vspace{5mm}

  \par

\vspace{5mm}

\section{Definitions. Notations. Statement of problem. Previous results.}

\vspace{5mm}

 \ Let $ \ D \ $ be certain proper own: $ D \ne R^n \ $  sub - domain of the whole space $ \ R^n, \ $ equipped with ordinary
 Euclidean norm $ \ |x|, \ x \in R^n, \ $
 for instance open, connected, convex, having non - trivial interior and Lipschitz boundary $ \ \partial D. \ $  The distance between arbitrary
 vector $ \ x \in R^n \ $  and boundary $ \ \partial D \ $  will be denoted as follows
$$
 d(x) = d(x,\partial D) = \inf_{y \in \partial D} |x - y|.
$$
 \ Let also $ \ u = u(x), \ x \in D \ $ be arbitrary valued function belonging to the Sobolev's space $ \ W^1(p): \ $

$$
||u||W^1_p \stackrel{def}{=} \left| \ \nabla u \right|_p^{1/p} + |u|_p < \infty,
$$
where the usually Lebesgue - Riesz norm $ \ |u|_p \ $ has the form

$$
|u|_p := \left[ \ \int_D |u(x)|^p \ dx   \ \right]^{1/p}, 1 \le p < \infty;
$$
$ \ \nabla(u) \ $ denotes the gradient in an distributional sense of the function $ \ u. \ $ \par

 \vspace{4mm}

 \ The following (and close)  key {\it weight}  Hardy - Sobolev - Poincare - Wirtinger (briefly, HSPW)
  inequality

\begin{equation} \label{HSPW}
\int_D \frac{|\nabla u|^p \ dx}{d^{\alpha}(x,\partial D)} \ge \left( \  \frac{\alpha + p - n}{p}   \ \right)^p \cdot
 \int_D \frac{|u|^p \ dx}{d^{p + \alpha}(x,\partial D)}, \ p > n - \alpha, \alpha = \const < n
\end{equation}
belongs to many authors, see e.g.  \cite{Avkhadiev1}, \cite{Avkhadiev2}, \cite{Barbatis}, \cite{Devyver}, \cite{Marcus}, \cite{Matskewich}, \cite{Pinchover}.
See also \cite{Hardy}, \cite{Maz'ya}, \cite{Mitrinovic}, \cite{Ostrovsky4}, \cite{Sobolev}, \cite{Talenty} etc.\par

 \ It  may be rewritten as follows. Introduce the following (Borelian) measure $ \ \mu_{\alpha} = \mu_{\alpha,D} \ $ on the measurable  subsets of
 the source set  $ \ D \ $

\begin{equation} \label{mu alpha}
\mu_{\alpha, D} (dx) \stackrel{def}{=} \frac{dx}{d^{\alpha}(x, \partial D) }
\end{equation}
and put

\begin{equation} \label{Koef K}
K(p) = K(p; \alpha,n) := \frac{p}{p + \alpha -n}, \ p > n - \alpha.
\end{equation}
 \ Define also the (linear) operator defined on all the functions having as a domain of definition the set $ \ D \setminus \partial D: \ $

\begin{equation} \label{T oper}
T[u](x) \stackrel{def}{=} \frac{u(x)}{d(x,\partial D)}.
\end{equation}

\ Then when $ \ p > n - \alpha \ $

\begin{equation} \label{imp ineq}
||T [u]||L_p(D,\mu_{\alpha}) \le K(p) \cdot ||\nabla u||L_p(D, \mu_{\alpha}).
\end{equation}

 \ Herewith the "constant" $ \ K = K(p) \ $ is the best possible for arbitrary proper domain $ \ D: \ $

\begin{equation} \label{best possible}
\sup_{u: \  ||\nabla u||L_p(D, \mu_{\alpha}) = 1 } ||T[u]||L_p(D,\mu_{\alpha})  = K(p), \ p > n - \alpha.
\end{equation}

 \vspace{4mm}

 \ {\bf  Our claim in this short preprint is to extrapolate the last estimate onto the so- called Grand Lebesgue Spaces
 instead the classical Lebesgue - Riesz ones, builded on the our set $ \ D \ $ equipped with the measure} $ \ \mu_{\alpha,D}. \ $ \par

 \vspace{4mm}

  \ {\it We intent herewith to calculate the exact value of correspondent embedding constants.}  \par

 \vspace{3mm}

 \ The particular case of these statement was considered in  \cite{Ostrovsky4}; see also
\cite{Ostrovsky1} - \cite{Ostrovsky3}. \par

\vspace{5mm}

 \ Let us recall here for readers convenience some known definitions and  facts  from  the theory of Grand Lebesgue Spaces (GLS)
adapted   exclusively to offered article.

   \ Let  the number $ \ b  \ $ be constants  such that $ \ n - \alpha  < b \le \infty; \ $   and let
   $ \ \psi = \psi(p) = \psi[b](p), \ p \in (n - \alpha,b), \ $ be numerical valued strictly positive function  not necessary to be finite in every point:

\begin{equation} \label{Positive psi}
 \inf_{p \in (n - \alpha,b)} \psi[b](p) > 0.
\end{equation}

\vspace{3mm}

 \ The set of all such a functions $ \ \psi(\cdot) \ $  will be denoted  by $ \ \Psi(b) = \{\  \psi = \psi(\cdot)  \ \}. \ $ \par

\vspace{3mm}

  \ For instance

$$
    \psi_m(p) := p^{1/m}, \ m = \const > 0,
$$
or more generally

\begin{equation} \label{more general}
\psi_{m,\beta,L}(p) := p^{1/m} \ \ln^{\beta}(p + 1) \ L(\ln ( p +1)), \ \beta = \const,
\end{equation}
where  $ \ L(\cdot) \ $ is positive continuous slowly varying as $ \ p \to \infty \ $ function.\par

\vspace{3mm}

 \ {\bf Definition 1.1.}

\vspace{3mm}

 \ By definition, the (Banach) Grand Lebesgue Space (GLS)    $  \ G \psi  = G\psi [b],  $
    consists on all the real (or complex) numerical valued measurable functions
   $   \  f: D  \to R \ $  defined on the whole our  space $ \ \Omega \ $ and having a finite norm

 \begin{equation} \label{norm psi}
 || \ f \ ||G\psi = ||f||G\psi[b] \stackrel{def}{=} \sup_{p \in (n - \alpha,b)} \left[ \frac{||f||_{p, D, \mu_{\alpha}}}{\psi(p)} \right].
 \end{equation}

 \vspace{4mm}

 \ The function $ \  \psi = \psi(p) = \psi[b](p) \  $ is named as  the {\it  generating function } for this space. \par

  \ If for instance

$$
  \psi(p) = \psi^{(r)}(p) = 1, \ p = r;  \  \psi^{(r)}(p) = +\infty,   \ p \ne r,
$$
 where $ \ r = \const \in [1,\infty),  \ C/\infty := 0, \ C \in R, \ $ (an extremal case), then the correspondent
 $ \  G\psi^{(r)}(p)  \  $ space coincides  with the classical Lebesgue - Riesz space $ \ L_r = L_r(\Omega, {\bf P}). \ $ \par

 \vspace{4mm}

  \ The finiteness of some GLS $ \ G\psi \ $ norm for the function  $ \ f:  D \to R \ $  is closely related with its tail function

$$
T_{f}(u) :=  \mu_{\alpha} \{x: |f(x)| > u \}, \ u \ge 0.
$$

 \ The GLS spaces are also closely related with the suitable exponential Orlicz ones, builded on the our measurable space
 $ \ (D,\mu_{\alpha}). \ $ \par

 \ See the detail investigation of these spaces in the works
  \cite{Ermakov}, \cite{Fiorenza1} - \cite{Fiorenza2}, \cite{fioforgogakoparakoNAtoappear}, \cite{fioformicarakodie}, \cite{Fiorenza-Formica},
  \cite{formicagiovamjom},  \cite{Kozachenko1} - \cite{Kozachenko3},  \cite{Ostrovsky0}, chapters 1,2.\par

\vspace{5mm}

\section{Main result.}

\vspace{5mm}

 \ {\bf Theorem 2.1.} Suppose that for certain function $ \ u = u(x), \ x \in D \ $  its gradient belongs to some Grand Lebesgue Space
 $ \ G \psi[b]: \ $

\begin{equation} \label{belongs Gpsi}
||\nabla u||G\psi[b] < \infty.
\end{equation}

 \ Introduce an auxiliary such a function

 $$
 \psi_K(p) := K(p) \cdot \psi(p), \ n-\alpha < p < b.
 $$

 \vspace{3mm}

 \ Our proposition:

\begin{equation} \label{main}
||T[u]||G\psi_K \le 1 \times ||\nabla u||G\psi,
\end{equation}
where  the constant "1" in (\ref{main}) is the best possible. \par

\vspace{3mm}

\ {\bf Remark 2.1.} The generating function $ \ \psi[b](\cdot) \ $  in (\ref{belongs Gpsi}) may be choosed by the
natural way:

$$
\psi[b](\cdot) := ||\nabla u||_p, \ n - \alpha < p < b,
$$
of course, if there exists such a value $ \ b > n - \alpha. \ $  \par

\vspace{4mm}

 \ {\bf Proof of upper bound} is very  simple and is alike to one in the articles
\cite{Ostrovsky2}, \cite{Ostrovsky3}; as well as the proof of the lower bound. Suppose

$$
||\nabla u ||G\psi < \infty;
$$
one can assume without loss of generality

$$
||\nabla u ||G\psi = 1.
$$

 \ It follows from the direct definition of the norm in GLS

$$
\forall p \in (n - \alpha, b) \ \Rightarrow  \ ||\nabla u ||_p \le \psi(p).
$$
 \ We apply the inequality (\ref{imp ineq})

$$
||T [u]||L_p(D,\mu_{\alpha}) \le K(p) \cdot \psi(p) = \psi_K(p),
$$
and on the other words

$$
||T[u]||G\psi_K \le 1 = ||\nabla u||G\psi.
$$

\vspace{3mm}

 {\bf The lower bound} in  (\ref{main}) follows immediately from one of the results of the article
\cite{Ostrovsky3}; see also \cite{Ostrovsky2}, taking into account the {\it exactness} of the value
$ \ K = K(p), \ $  see (\ref{best possible}). \par
Q.E.D. \par

\vspace{5mm}

\section{Necessary conditions for these estimations.}

\vspace{5mm}

  \ We will ground in this section that the  "configuration" given by the estimate (\ref{HSPW}) is essentially
  non - improvable. Namely, suppose that there exists non - trivial constants $ \ a,h; G(D,a,h) \ $ such that
  for {\it arbitrary} proper non - trivial sub - domain $ \ D \subset R^n \ $
 and for  {\it arbitrary} function $ \ u = u(x), \ x \in D \ $ belonging  to the space $ \  C^0_{\infty}(D) \ $ there holds
 the estimate

\begin{equation} \label{pq ineq L}
 L[u] := L[u](a) \stackrel{def}{=}    \left[ \ \int_D \frac{|u(x)|^q}{d^a(x,\partial D)} \ dx \  \right]^{1/q} \le
\end{equation}

\begin{equation} \label{pq ineq R}
R[u]:= R[u](h) \stackrel{def}{=} G(D,a,h) \ \left[ \ \int_D \frac{|\nabla u|^p}{d^h(x,\partial D)}  \ dx \ \right]^{1/p}
\end{equation}
for certain fixed values $ \ p,q \in [1,\infty). \ $ \par

\vspace{4mm}

 \ {\bf  Theorem 3.1. } Suppose that the relations   (\ref{pq ineq L}) and (\ref{pq ineq R}) holds true for arbitrary proper domain $ \ D \ $
  as well as for certain  non - zero  function  $ \ u = u(x), \ x \in D \ $ from the set $ \ C^0_{\infty}(D). \ $  Then

\begin{equation} \label{key relation}
\frac{a - n}{q} = 1 + \frac{h-n}{p}.
\end{equation}

\vspace{4mm}

 \hspace{3mm} {\bf Proof.} We will apply the so - called dilation method, belonging at first, perhaps, to
G.Talenty, see the article \cite{Talenty}, in which was considered the particular case $ \ a = h = 0. \ $ \par
 \ Let us choose the domain $ \ D \ $ on the form

$$
D := R^d \otimes R^r_+,  \ n = \dim D = d + r, \  r,d > 0.
$$

 \ Let's agree to write instead of $ \ x \ $ the vector $ \ (x,y), \ x \in R^d, \ y \in R^r_+. \ $ Then

$$
L^q[u]  = \int_{R^d} \int_{R^r_+} \frac{|u(x,y)|^q \ dx \ dy}{|y|^a},
$$

$$
R^p[u] =  \int_{R^d} \int_{R^r_+} \frac{|\nabla(x,y)|^p \ dx \ dy}{|y|^h}.
$$
 \ One can choose the function $ \ u = u(x,y) \ $ from the set $ \ C^0_{\infty}(D)  \ $ such that
 $ \ L(u) > 0, \ R[u] > 0;  \  $ therefore $ \ 0 <L[u] \le C R[u] < \infty. \ $  \par

  \ Let also $ \ \lambda \in (0,\infty) \ $ be a parameter of dilatation:

$$
V_{\lambda}[u](x,y) \stackrel{def}{=} u(\lambda x, \lambda y).
$$
  \ Obviously, $ \ V_{\lambda}[u](\cdot, \cdot) \in C^0_{\infty}(D). \ $ Therefore, it satisfies the inequality
(\ref{pq ineq L}) \ - \  (\ref{pq ineq R}):

\begin{equation} \label{lambda}
L_a[V_{\lambda} u] \le C R_h[V_{\lambda} u].
\end{equation}

 \ But

$$
\left\{ L_a[V_{\lambda} u] \right\}^q = \lambda^{-d - r + a} L^q[u];
$$

$$
 L_a[V_{\lambda} u] = \lambda^{(a - d - r)/q } L_a[u]
$$

and analogously

$$
R_h[V_{\lambda} u] = \lambda ^{(p - d - r + h)/p} R_h[u],
$$
and we get to the inequality

$$
C_1 \lambda^{(a - d - r)/q } \le C_2 \lambda ^{(p - d - r + h)/p}
$$
for arbitrary positive values $ \ \lambda. \ $  Following,

$$
\frac{ a - d - r}{q}   =  \frac{p - d - r + h}{p},
$$
which is equal to (\ref{key relation}). \par

 \ Of course, in the case when $ \ a\cdot h = 0 \ $  or  $ \ p = q \ $ we get to the classical Sobolev's
 inequality, see e.g.  \cite{Talenty}. \par

\vspace{5mm}

\section{Concluding remarks.}

\vspace{5mm}

 \hspace{3mm} {\it Open question: } Suppose that the there holds (\ref{key relation}) and that the
domain $ \ D \ $ is proper. Assume also the function $ \ u = u(x) \ $ belongs to the set $ \ C^0_{\infty}(D). \ $
Will be the inequality (\ref{pq ineq L}) - (\ref{pq ineq R}) true, of course,
in general case, i.e.  when $ \ p \ne q? \ $ \par

\vspace{6mm}

\vspace{0.5cm} \emph{Acknowledgement.} {\footnotesize The first
author has been partially supported by the Gruppo Nazionale per
l'Analisi Matematica, la Probabilit\`a e le loro Applicazioni
(GNAMPA) of the Istituto Nazionale di Alta Matematica (INdAM) and by
Universit\`a degli Studi di Napoli Parthenope through the project
\lq\lq sostegno alla Ricerca individuale\rq\rq (triennio 2015 - 2017)}.\par

\end{document}